\documentclass[12pt]{article}

\def\@fnsymbol#1{\ifcase#1\or \@arabic\c@footnote\or \@arabic\c@footnote\or \@arabic\c@footnote\or \@arabic\c@footnote\or \@arabic\c@footnote\else\@ctrerr\fi}
\makeatother

\usepackage[english]{babel}
\usepackage[utf8]{inputenc}
\usepackage[T1]{fontenc}
\usepackage[top=90pt,bottom=80pt,left=75pt,right=75pt]{geometry}
\usepackage{mathtools, amsmath, amsthm, amsfonts, amssymb}
\usepackage{array}
\usepackage{nicematrix}
\usepackage{graphicx}
\usepackage{xargs}
\usepackage{enumerate}
\usepackage{csquotes}
\usepackage{fancyhdr}
\usepackage{xurl}
\usepackage{hyperref}
\usepackage[backend=biber, style=numeric, giveninits=true, 
 doi=false, isbn=false, url=false, hyperref=true, maxbibnames=99, minbibnames=99]{biblatex}
 
\addtocounter{footnote}{2}

\newtheorem{example}{Example}

\newcommand{\R} {\mathbb{R}}
\newcommand{\bR} {\mathbb{R}}

\newtheorem{theorem}{Theorem}[section]

\newtheorem{proposition}[theorem]{Proposition}
\theoremstyle{definition}

\theoremstyle{definition}
\newtheorem{definition}[theorem]{Definition}
\theoremstyle{remark}

\newcommand{\beqo}{\begin{eqnarray*}}
\newcommand{\beq}{\begin{eqnarray}}
\newcommand{\eeqo}{\end{eqnarray*}}
\newcommand{\eeq}{\end{eqnarray}}

\newcommand{\mat}{\left[ \begin{array}{c}  }
\newcommand{\rix}{\end{array} \right] }

\numberwithin{equation}{section}

\newcommand{\lw}[1]{\smash{\lower2.ex\hbox{#1}}}

\begin{document}
\title{Factorized sparse approximate inverse preconditioning for singular M-matrices
}
\author{Katherina Bick\footnote{katherina.bick@gmx.de} \and  Reinhard Nabben\thanks{Corresponding author, R. Nabben, TU Berlin, Institut f\"ur Mathematik, MA 3-3, Stra{\ss}e des 17. Juni 136, D-10623 Berlin, Germany, phone: +49-30-314 29291, (nabben@math.tu-berlin.de)}} 
\date{}

\maketitle

\begin{abstract}
Here  we consider the factorized sparse approximate inverse (FSAI) preconditioner introduced by Kolotilina and Yeremin in \cite{KolY93}. We apply the FSAI preconditioner to singular irreducible  M-matrices. These  matrices arise e.g. in discrete Markov  chain modeling or as graph Laplacians.  We show, that  there are some restrictions on the nonzero pattern needed for a stable  construction of the FSAI preconditioner in this case. 
With these restrictions  FSAI  is well-defined. Moreover, we proved that the FSAI preconditioner shares some important properties with the original system. The lower triangular matrix $L_G$ and the upper triangular matrix $U_G$, generated  by FSAI,  are non-singular and non-negative. The diagonal entries of $L_GAU_G$  are positive and  $L_GAU_G$,  the preconditioned matrix, is a singular M-matrix. Even more, we establish  that  a (1,2)-inverse is computed  for the complete nonzero patter.
\end{abstract}


{\bf AMS:} 
65F10, 65F50, 65N22, 65N55 \\

\begin{center}
 {\small \it Dedicated to Daniel B. Szyld on his 70th birthday}
\end{center}

\pagestyle{myheadings}
\thispagestyle{plain}
\markboth{K. BICK AND R. NABBEN}{FSAI for singular M-matrices}

\section{Introduction}
We consider  the iterative solution of linear systems of the form
\begin{eqnarray}
Ax = b  \quad A \in \mathbb{C}^{n \times n}, \quad x,b \in \mathbb{C}^{n}, \label{lineq}
\end{eqnarray}
where $A$ is a  large, sparse, and possibly non-symmetric matrix. 
Krylov subspace methods with preconditioning are nowadays  the methods of choice to solve \eqref{lineq}. The first preconditioners like the Jacobi-method or the incomplete LU-factorization are invented in the late seventies. These methods  are so-called implicit preconditioners. The preconditioner $M$ approximates the matrix $A$
and inside the Krylov subspace method a linear system with $M$ has to be solved in each iteration.
In the nineties so-called explicit   preconditioner became  popular. 
Here the preconditioner $B$ approximates the inverse of $A$ and thus just a multiplication with the preconditioner $B$ is needed in the Krylov subspace method. 
Since $B$  should be simple to compute and should be sparse, a nonzero pattern for $B$ is often used. This nonzero pattern can be chosen in advance or can be created dynamically during the computation of $B$.  Among  these  approximate inverse preconditioner are the SPAI, the FSAI, the FSPAI  and the  AINV method \cite{KolY93,BenMT96,GroH97,BenT98a,BenT98,BenCT00, KhaKNY01, BenT03}. These sparse approximate inverse preconditioners also have natural parallelism which  make them even more attractive. 
For review papers on preconditioning   we refer to \cite{BenT99, Ben02}. 

Approximate inverses are also used for shifted systems, as smoother in multigrid methods,  and in regularization and parallel computing, see e.g. \cite{ BarBS99, BenB03, BolHK16, BroGMR01, HucS10, Huc00, Huc03, AnzHBD18}.


Here we consider singular matrices. Theoretical results  for  different iterative methods  for singular matrices can be found e.g. in \cite{Szy94, MarS00, NABS07,FroNS08,LudNT16,KehN16,ErlN17,ElsFNSS03}. Moreover we concentrate on  irreducible singular M-matrices.
Non-negative matrices and (singular) M-matrices and their
generalized or group inverses  are useful tools not only in matrix analysis, but also in the analysis of several applications. Among them  are   stochastic processes e.g. Markov models, graph theory - graph Laplacian, electrical networks, and demographic models see e.g. \cite{BerP79,Mey82,CamM09,KirN13}. Discrete Markov chains arise in reliability modeling, queuing network analysis, large scale economic modeling and computer system performance evaluation \cite{BenT02}. For different iterative methods  for Markov chains and M-matrices see e.g. \cite{BruPS05, MarS93, BenNFS01, MenN08a, MenN08b, ElsFNSS03}. For  some algebraic multigrid methods for graph  Laplacian (i.e. singular M-matrices) see e.g. \cite{BolFFHK11, LivB12, NapN16}. 

The application of approximate inverse techniques to singular systems raises several interesting questions. First, are the algorithms stable for all nonzero pattern? Second, since  the inverse of $A$ does not exists, it is not clear what matrix is then approximated.



The SPAI  and the AINV preconditioner   work for  irreducible singular M-matrices. For the complete nonzero pattern SPAI leads  to the  Moore-Penrose inverse of $A$ while  AINV gives a (1,2)-inverse, see e.g. \cite{BenT02}. 
But the factorized sparse approximate inverse preconditioner (FSAI) introduced by Kolotilina and Yeremin in \cite{KolY93} 
(see also \cite{KolY95, YerK00, KolY99}) 
is not yet considered for non-singular M-matrices.

Kolotilina and Yeremin showed in \cite{KolY93} that FSAI works for Hermitian positive definite matrices and non-singular M-matrices and non-singular H-matrices, but singular matrices  are not considered. 

In \cite{FerJP14})  the authors  generalize to the unsymmetric case the Block
Factorized Sparse Approximate Inverse (Block FSAI) studied so far for Hermitian positive definite  matrices, see \cite{JanFG10, JanF11}. In their algorithm systems  with  singular matrices can occur, which leads to a breakdown of the algorithm. In this case, the solution to the corresponding system is skipped setting some relevant vectors to the null vector (see page 235 in \cite{FerJP14}).

In this  note  we consider  the FSAI for singular irreducible  M-matrices. We show, that  there are some restrictions on the nonzero pattern needed for a stable  construction of the FSAI preconditioner in this case. 
With these restrictions  FSAI  is well-defined and does not break down. So we closed the gap, i.e. not only SPAI and  AINV but also  FSAI can be applied to singular irreducible M-matrices without any breakdown.

Moreover,  FSAI preconditioner shares some important properties with the original system. We prove that the lower triangular matrix $L_G$ and the upper triangular matrix $U_G$, generated  by FSAI,  are non-singular and non-negative. The diagonal entries of $L_GAU_G$  are positive and  the preconditioned matrix $L_GAU_G$  is a singular M-matrix. Even more, we establish  that  a (1,2)-inverse is computed  for the complete nonzero patter.

This note is organized as follows.  In section 2 we list some properties  of irreducible M-matrices. Section 3 then briefly describes the FSAI preconditioner. The main results are then given in section 4. 
In section 5  we give  some small  numerical examples. 

\section{Singular M-matrices and generalized inverses}

The {\it Moore-Penrose inverse} of 
$A \in \mathbb R^{n,n}$ is a  matrix $A^{+}$ which satisfies:
\begin{center}$\begin{array}{lrcl}
(1)& AA^{+}A&=&A\cr
(2)& A^{+}AA^{+}&=&A^{+}\cr
(3)& (AA^{+})^T&=&AA^{+}\cr
(4)& (A^{+}A)^T&=&A^{+}A.\cr
\end{array}$ \end{center}

Every  matrix $A\in \mathbb R^{n,n}$ has a unique Moore-Penrose inverse.
If a matrix $A^{\dagger}$ satisfies conditions (1) and (2), then the matrix is called an (1,2)-inverse. This kind of generalized inverse is not unique. Indeed there may infinitely many such inverses. For more details on generalized inverses see \cite{CamM09,KirN13} and references therein.

Next  we list  some properties of (singular) M-matrices which  we need in this note. First we give   the definition of Z-matrices. 

\begin{definition} A real matrix $A = [a_{i,j}]$ is called Z-matrix, if $a_{i,j} \leq 0$ for $i \neq j$.
 \end{definition}
 
The class  of Z-matrices was  defined by Fiedler and Ptak in \cite{FieP62}. Later in \cite{FieM92} a more detailed analyses  of different  classes of Z-matrices is given, see also \cite{NabV95,Nab97}.  
 
 Among  the classes  of Z-matrices, the most prominent class is the  class of (non-singular) M-matrices, since M-matrices arise in many applications, see e.g. \cite{BerP79}. 
 
 \begin{definition} A Z-matrix $A = [a_{i,j}]$ is called (non-singular) M-matrix, 
 if  $A$  can  be written  as 
 \begin{equation*}
  A = sI - B, 
 \end{equation*}
where $I$ is the identity matrix, $B$ is a (entry wise) non-negative matrix and $s \in \R$ with 
 \begin{equation*}
s > \rho(B).
 \end{equation*}
 \end{definition}
 
 In \cite{BerP79}  there  are more than 50 equivalent conditions given for a Z-matrix to be an M-matrix, see  Theorem 2.3 in \cite{BerP79}.  We use the following ones:
 
 \begin{theorem} \label{theo:M-ma}
Let $A\in R^{n,n}$ be a  non-singular Z-matrix. Then  the following are equivalent:
\begin{enumerate}
\item[(1)] $A$ is an M-matrix, 
\item[(2)] There exists an  entry wise positive vector $x$, i.e. $x>0$, such that  $Ax$ is entry wise positive, i.e. $Ax > 0$
\item[(3)] $A^{-1}$ is entry wise non-negative, i.e. $A^{-1} \geq 0$. 
\end{enumerate}
\end{theorem}
 
 \begin{definition} A Z-matrix $A = [a_{i,j}]$ is called singular M-matrix, 
 if  $A$  can  be written  as 
 \begin{equation*}
   A = sI - B, \  \text{with} \ B \  \text{non-negative and} \  s =  \rho(B).
 \end{equation*}
 \end{definition}
 
 There  are also many equivalent conditions known for a singular Z-matrix   to be singular M-matrix (see Theorem 4.6  in \cite{BerP79}). We just need the following one.

 \begin{theorem}\label{theo:singM-ma}
Let $A\in R^{n,n}$ be a  singular Z-matrix. Then  the following are equivalent:
\begin{enumerate}
\item[(1)] $A$ is a  singular M-matrix, 
\item[(2)]  For every $\epsilon > 0$, the matrix $A + \epsilon I$ is a non-singular M-matrix.
\end{enumerate}
\end{theorem}
 
 Here  we  consider irreducible singular M-matrices. We then have some more properties, see Theorem 4.16 in\cite{BerP79}.
\begin{theorem}\label{theo:irresingM-ma}
Let $A\in R^{n,n}$ be a singular irreducible M-matrix. Then:
\begin{enumerate}
\item[(1)] $A$ has rank $n-1$, 
\item[(2)] There exists a  vector $x>0$ with $Ax=0$. 
\item[(3)] Every principal submatrix of order at most $n-1$ is a non-singular M-matrix.
\item[(4)] There exist a  lower and a upper triangular matrix  $L_A$ and $ U_A\in R^{n,n}$ with ones one the diagonal and a diagonal matrix 
$D=\,$diag$\,(d_1,\ldots,d_{n-1},0)$ with $d_1,\ldots,d_{n-1}$ positive, such that $A=L_ADU_A$.
\end{enumerate}
\end{theorem}

\section{FSAI preconditioning}
Here  we give a short introduction into the {FSAI} procedure introduced in \cite{KolY93}  for non-singular matrices. 
Note, that   we prescribe nonzero pattern $S_L$ and $S_U$.  

First choose  arbitrary nonzero pattern
\begin{eqnarray*}
\hat S_L\subseteq\{(i,j):i\ge j\}, \quad \quad \hat S_U\subseteq\{(i,j):i\le j\}.
\end{eqnarray*}

The only condition we need  to have is that the pairs $(i,i)$ must be in these nonzero pattern.  So we define  
\[S_L = \hat S_L  \cup \{(i,i) : i = 1, \ldots , n \},  \quad \quad  S_U = \hat S_U  \cup \{(i,i) : i = 1, \ldots , n \}. \]

Then compute the entries of the lower triangular matrix $L_G$ and the entries of the upper triangular matrix $U_G$  such that
\[(L_GA)_{ij}=\delta_{ij}, \ \text{for} \quad (i,j) \in S_L; \quad \quad (AU_G)_{ij}=\delta_{ij}, \quad \text{for} \ (i,j) \in S_U. \]
The entries at positions $(i,j) \notin S_L$ or  $(i,j) \notin S_U$ in $L_G$ and $U_G$ are zero. 

These  conditions lead to linear systems for the entries of $L_G$ and $U_G$, which  can be solved in parallel. Once the matrices $L_G$ and $U_G$ are computed, 
one builds the diagonal matrix $D$ with   $D=\,$diag$\,(L_GAU_G).$

If $D$ has positive entries on the  diagonal we obtain the left and right preconditioner
\[G_1:=D^{-\frac{1}{2}}L_G \quad \mbox{and} \quad G_2:=U_G D^{-\frac{1}{2}}.\]

It is proved in \cite{KolY93} that for non-singular M-matrices the FSAI procedure is well-defined for all nonzero patterns $S_L$ and $S_U$. Moreover, $L_GAL_U$ is a non-singular M-matrix. 

So  we can ask the following questions. Is the FSAI procedure also well-defined for singular M-matrices or are there  some modifications needed?
Moreover, if the FSAI procedure works, which matrix, i.e. which generalized  inverse will be approximated?

The  first question has a negative answer. For singular matrices some modifications are needed. 

Find  $L_G$ with  $(L_GA)_{ij}=\delta_{ij} \ \ {\mbox{for}} \ \
(i,j) \in S\subseteq\{(i,j):i\ge j\} $
is generally not possible which is shown by the following counter example. 

\begin{example} \label{ex:1}
 Let $A= \left[ \begin{array}{rr}3 & 
-3 \\  -3 & 3 \end{array} \right]$. $A$ is a  singular, irreducible 
M-matrix. We  consider the complete nonzero pattern
\begin{center}
$S =\{(i,j):i \ge j\}$\\
\end{center}

So we  want: $L_G=\left[ \begin{array}{rr}l_{11} & 0\cr l_{21} & l_{22}\end{array}\right]$ with 
$L_GA = \left[ \begin{array}{rr}1 & \ast\cr 0 & 1\end{array}\right]$.
But this leads to the equations
\begin{eqnarray*} 3 l_{21} -3 l_{22} & = & 0\\
  -3  l_{21} +3 l_{22} & = & 1
\end{eqnarray*}\\
which can not be solved.
\end{example}

Example \ref{ex:1}  shows  that the FSAI algorithm does not work  for singular matrices with  complete  nonzero patterns. But this is not a surprise. In general we have:

Let $a_i$ be the columns of $A$ and let $l_j^T$ be the rows of $L$. If $A$ is a singular irreducible M-matrix,  we have
\[
 a_n = \sum_{i=1}^{n-1} \lambda_ia_i
\]
for some $\lambda_i \in \bR$, (Theorem \ref{theo:irresingM-ma}, (2)). 
If we choose the complete nonzero pattern then  we need to have
\[
 l_n^Ta_i = 0 \quad \text{for} \quad i = 1, \ldots, n-1.
\]
and 
\[
 l_n^Ta_n = 1.
 \]
 But this contradicts
 \[
 l_n^Ta_n =  l_n^T \sum_{i=1}^{n-1} \lambda_ia_i = \sum_{i=1}^{n-1} \lambda_i l_n^Ta_i = 0.
 \]

\section{Main results}
First we will  show, that some simple restrictions on the nonzero pattern lead to a stable FSAI method for singular irreducible M-matrices. 

Again, we  start  with some arbitrary subsets of pairs of indices. But we exclude  the pairs $(n,n-1)$ and $(n-1,n)$.
\[\tilde S_L\subseteq\{\{(i,j):i\ge j\}\,\backslash \,(n,n-1)\}, \quad \quad 
\tilde S_U\subseteq\{\{(i,j):i\le j\}\,\backslash \,(n-1,n)\}\]

then, as above,  we have to include the diagonal entries
\begin{equation} \label{eq:pattern}
 S_L = \tilde S_L  \cup \{(i,i) : i = 1, \ldots , n \},  \quad \quad S_U = \tilde S_U  \cup \{(i,i) : i = 1, \ldots , n \} 
\end{equation}
which gives our nonzero pattern. 

We then want to find a lower triangular matrix  $L_G$ and an upper triangular matrix $U_G$ with
\begin{eqnarray} \label{fspai1}
(L_GA)_{ij} & = & \delta_{ij} \ \ {\mbox{for}} \ \ (i,j) \in S_L \\ \label{fspai2}
(AU_G)_{ij} & = & \delta_{ij} \ \ {\mbox{for}} \ \ (i,j) \in S_U.
\end{eqnarray}
If $(i,j) \notin S_L$ then $(L_G)_{ij} = 0$  and if   $(i,j) \notin S_U$ then $(U_G)_{ij} = 0$.
Define 
  \begin{equation} \label{fspaiD}
 D := diag(d_1, \ldots , d_{n-1}, d_{n,n}) :=  \ diag(L_GAU_G)
\end{equation}

and as above build  preconditioners $D^{-\frac{1}{2}}L_G,$ and  $U_G D^{-\frac{1}{2}}$, if the diagonal entries of $D$ are positive.

We   will show  that now the FSAI  procedure is well-defined and leads to a singular M-matrix. Moreover, we will  prove  that  for  the complete  nonzero pattern with 
\begin{equation*}
 D^- = diag(d_1^{-1}, \ldots , d_{n-1}^{-1}, 0),
\end{equation*}

the matrix $U_GD^-L_G$ is a (1,2)-inverse of $A$. 

We  will also see that  any other pair $(n,i)$  and $(i,n)$  with $i \neq n$  instead of  $(n,n-1)$ and $(n-1,n)$ can  be excluded in the nonzero pattern. Note that this  restriction is not a drawback in praxis, since  mostly sparse matrices will be constructed. But for the theoretical results   we need to keep  that in mind.    

We  start  with the following proposition.
\begin{proposition}
 Let $A$  be a singular irreducible M-matrix. Then there exists for all nonzero pattern $S_L$ and $S_U$ as in \eqref{eq:pattern} a unique lower triangular matrix $L_G$ that satisfies \eqref{fspai1} and a unique upper triangular matrix $U_G$ that satisfies \eqref{fspai2}.
\end{proposition}
\begin{proof}
 The rows of $L_G$ can be computed in parallel.  For the entries of each row we need to solve a linear system. But the coefficient matrices of these systems  are principal submatrices of $A$  with maximal order $n-1$. But these principal submatrices of $A$ are non-singular M-matrices by Theorem \ref{theo:irresingM-ma}. Thus  each linear system has a unique  solution.  
 Hence, $L_G$ is well-defined and unique. The same  holds for $U_G$. 
 \end{proof}

\begin{theorem} \label{theo:exis}
 Let $A$  be a singular irreducible M-matrix. Let $L_G$ and $U_G$ be the  unique matrices given by  \eqref{fspai1} and \eqref{fspai2}. Then 
$L_G$ and $U_G$ are (entry-wise) non-negative and non-singular.  
 \end{theorem}
\begin{proof}
 Let $B \in \bR^{n-1,n-1}$ be the principal submatrix of $A$ obtained from $A$  by deleting the last  row and column of $A$. Since $A$ is an irreducible M-matrix, $B$ is a non-singular M-matrix, see  Theorem \ref{theo:irresingM-ma}.
 
 Define $B^{(i)}(S_L)\in \mathbb R^{n-1,n-1}$ for $i=1,\ldots,n-1$ by \\
 \[
 (B^{(i)}(S_L))_{rs}=  
\left\{ \begin{array}{cl}  a_{rr}, & r=s \\
                           a_{rs}, & (i,r) \in S_L \wedge (i,s) \in S_L \\
                              0, & \mbox{else.}\cr
        \end{array}
 \right.\]

Since $B$ is a non-singular M-matrix  the matrices $B^{(i)}(S_L)$ are also  non-singular M-matrices for all $i$. Hence, $(B^{(i)}(S_L))^{-1}$ exist for $i=1,\ldots,n-1$ and $(B^{(i)}(S_L))^{-1} \geq 0$.
Then define $L_F$ by  $(L_FB)_{ij}=\delta_{ij}$ for $(i,j) \in S_L\,\backslash\,\{(n,j):j=1,\ldots,n\}.$

But  the nonzero entries of the $i$-th row of $L_F$ equal the corresponding entries of the $i$-th row of $(B^{(i)}(S_L))^{-1}$. Thus the first $n-1$ rows of $L_G$ are non-negative. 
If we choose  $B \in \bR^{n-1,n-1}$ as the principal submatrix of $A$ obtained by deleting the first row and column of $A$, then we obtain that last row of $L_G$ is also non-negative. 
Since  the matrices $B^{(i)}(S_L)$ are M-matrices, the diagonal entries of $L_G$ are positive. Thus $L_G$ is a non-negative non-singular matrix. Similarly we prove that $U_G$ is a non-negative and non-singular matrix.
\end{proof}

\begin{theorem} \label{theo:neusingM}
 Let $A$  be a singular irreducible M-matrix. Let $L_G$ and $U_G$ be the  unique matrices given by  \eqref{fspai1} and \eqref{fspai2}. Then 
$L_GAL_U$ is a singular M-matrix.  
 \end{theorem}
 \begin{proof}
  With Theorem  \ref{theo:irresingM-ma}  there exists a positive vector $v$ such that $Av = 0$. 
  Thus we  obtain  $L_GAv =  0$.
  Hence
  \begin{equation*}
   (L_GA + \epsilon I)v = L_GAv + \epsilon v > 0 \quad \text{for all} \  \epsilon > 0.
  \end{equation*}
Now consider the off-diagonal entries of $(L_GA + \epsilon I)$ for an arbitrary $\epsilon > 0.$
We have 
\begin{equation*}
(L_GA + \epsilon I)_{i,j} = (L_GA)_{ij}  \quad \text{for} \quad i \neq j.
\end{equation*}
For $(i,j) \in S_L$, $i \neq j$ we have $(L_GA)_{ij} = 0$ by the construction of $L_G$. But for  $(i,j) \notin S_L$  we have 
 \begin{equation*}
 (L_GA)_{ij} = \sum_{k \leq i} (L_G)_{ik}(A)_{kj} \leq 0.
 \end{equation*}  
 The last inequality holds since in the  sum only the indices  $k$  are used for which $(i,k) 
\in S_L$. Otherwise the entries of $L_G$ are zero by construction. But we considered the case $(i,j) \notin S_L$, so $k \neq j$. Hence, since $A$ is a singular $M$-matrix,  the only positive entry $(A)_{kk}$ does not appear in the sum. 

Hence $L_GA + \epsilon I$ is a non-singular M-matrix for all $\epsilon > 0$. Thus $L_GA$ is a singular M-matrix, see Theorem \ref{theo:singM-ma}.

Similarly  one can prove that $AU_G$ is a singular M-matrix.

Next we consider  $(L_GAU_G + \epsilon I)_{ij}$  for $ i \neq j$. We then have 
$(L_GAU_G + \epsilon I)_{ij}$ = $(L_GAU_G)_{ij}$. For $i < j$ it holds
\begin{equation*}
 (L_GAU_G)_{ij} = \sum_{k=1}^{n} (L_G)_{ik}(AU)_{kj} = \sum_{i \geq k} (L_G)_{ik}(AU)_{kj} \leq 0,
\end{equation*}
since with $i < j$ and $i \geq k$ then $k \neq j$. 
Similarly
\begin{equation*}
(L_GAU_G)_{ij} = \sum_{k=1}^{n} (L_GA)_{ik}(U)_{kj} \leq 0.
\end{equation*}
Hence $L_GAU_G + \epsilon I$ is a non-singular M-matrix for all $\epsilon > 0$. Thus $L_GAU_G$ is a singular M-matrix. 
 \end{proof}

 So we have seen that $L_GAL_U$ is a singular M-matrix.  But  singular  M-matrices can have zeros on the diagonal. The next theorem excludes this.  
 
 \begin{theorem} \label{theo:nonsing}
 Let $A$  be a singular irreducible M-matrix. Let $L_G$ and $U_G$ be the  unique matrices given by  \eqref{fspai1} and \eqref{fspai2}. Then the  diagonal entries of 
$L_GAU_G$ are positive.  
 \end{theorem}
 \begin{proof}
As  in the proof of Theorem \ref{theo:exis}  let $B \in \bR^{n-1,n-1}$ be the principal submatrix of $A$ obtained from $A$  by deleting the last  row and column of $A$. Since $A$ is an irreducible M-matrix, $B$ is a non-singular M-matrix.  Moreover define $L_F$ and $U_F \in \bR^{n-1,n-1}$  by
\begin{equation*}
 (L_FB)_{kj} = \delta_{kj}  \quad \text{for} \quad (k,j) \in S_L \backslash \{(n,j) : j = 1, \ldots , n\}
\end{equation*}
and by 
\begin{equation*}
 (BU_F)_{kj} = \delta_{kj}  \quad \text{for} \quad (k,j) \in S_U \backslash \{(n,j) : j = 1, \ldots , n\}.
\end{equation*}
Next let $M$ be the leading principal submatrix of order $n-1$ of $L_GAU_G$. Then for $i \neq n$ and $j \neq n$
\begin{eqnarray*}
 (M)_{ij}  & =  & \sum_{k = 1}^{n} (L_G)_{ik}(AU_G)_{kj} \\
 & = &  \sum_{k = 1}^{n-1} (L_G)_{ik}(AU_G)_{kj} \\
 & = &  \sum_{k = 1}^{n-1} (L_G)_{ik}\left(\sum_{l=1}^{n}(A)_{kl}(U_G)_{lj}\right) \\
  & = &  \sum_{k = 1}^{n-1} (L_G)_{ik}\left(\sum_{l=1}^{n-1}(A)_{kl}(U_G)_{lj}\right) \\
  & = &  \sum_{k = 1}^{n-1} (L_G)_{ik}\left(\sum_{l=1}^{n-1}(B)_{kl}(U_G)_{lj}\right) \\ 
  & = & (L_FBU_F)_{ij}.
  \end{eqnarray*}
But now $B$ is a non-singular M-matrix, hence in this case we know from \cite{KolY93} that $(L_FBU_F)$  is also an M-matrix. Hence the first $n-1$ diagonal entries  of  $L_GAU_G$ are positive.

Similarly to $B$   we can use   the submatrix $F$ of $A$ obtained by deleting the $n-1$-th  row and column  and use the same  techniques  as above. This can be done since $(L_G)_{n,n-1} = (U_G)_{n-1,n} = 0$ which in turn allows to interchange
the last two rows and columns. We then have that the last  diagonal entry is positive also.
\end{proof}

Theorem \ref{theo:nonsing} guarantees that the diagonal  entries  of  $L_GAU_G$ are positive. Hence, with $D = diag(L_GAU_G)$, the left and right preconditioner $D^{\frac{1}{2}}L_G$ and  $U_GD^{\frac{1}{2}}$ are well  defined and non-singular.

In the  following we will  consider some properties of the above constructed matrices if the complete nonzero pattern is used, i.e. 
\begin{eqnarray} \label{S_L-voll}
S_L = \{(i,j):i\ge j\}\,\backslash \,\{(n,n-1)\}, \\ \label{S_U-voll}
S_U = \{(i,j):i\le j\}\,\backslash \,\{(n-1,n)\}. 
\end{eqnarray}

\begin{theorem} \label{theo:structure}
 Let $A$  be a singular irreducible M-matrix. Let $L_G$ and $U_G$ be the  unique matrices given by  \eqref{fspai1} and \eqref{fspai2} using $S_L$ and $S_U$ as given in \eqref{S_L-voll} and \eqref{S_U-voll}. Then the  matrix $L_GAU_G$ has the following form
 \begin{equation} \label{eq:IF}
   I_F := L_GAU_G = \left[ \begin{array}{ccccc}
                           d_1 & & & & \\
                           & \ddots & & & \\
                           & & d_{n-2} & & \\
                           & &  & d_{n-1} & b\\
                           & & & c & d_n
                          \end{array} \right],
 \end{equation}
 
where  the $d_1, \ldots , d_n$  are positive and $d_{n-1}d_n = cb.$
 \end{theorem}
\begin{proof}
Consider  the entries of $L_GAU_G$. For $i < j$ and $i < n - 1$ we have 
\begin{equation*}
 (L_GAU_G)_{ij} = \sum_{k=1}^n (L_G)_{ik}(AU_G)_{kj} = \sum_{ k \leq i} (L_G)_{ik}(AU_G)_{kj} = 0.
\end{equation*}
Similarly  we obtain  for  $i > j$ and $j < n - 1$
\begin{equation*}
 (L_GAU_G)_{ij} = \sum_{k=1}^n (L_GA)_{ik}(U_G)_{kj} = \sum_{ k \leq j} (L_GA)_{ik}(U_G)_{kj} = 0.
\end{equation*}
Thus $L_GAU_G$  has the  structure  as given in \eqref{eq:IF}.
 With Theorem \ref{theo:nonsing} the diagonal entries of $L_GAU_G$ are positive, thus $d_1, \ldots , d_n$  are positive. But  $L_GAU_G$ is a singular M-matrix (see Theorem \ref{theo:neusingM}) hence $d_{n-1}d_n = cb,$  which  completes the proof. 
\end{proof}

As mentioned in the beginning, excluding the pair $(n,n-1)$ and $(n-1,n)$ guarantees that the construction of FSAI is well-defined.  But  we can exclude any other pair $(n,j)$  and $(j,n)$  for $j \neq n$. Then  the matrix $I_F$ in \eqref{eq:IF} would have nonzero entries in the $(n,j)$  and $(j,n)$ positions  rather then in  $(n,n-1)$ and $(n-1,n)$. \\

Next define the matrix $D^- \in \bR^{n,n}$ by
\begin{equation}
 D^- = diag(d_1^{-1}, \ldots , d_{n-1}^{-1}, 0),
\end{equation}
where the $d_i$ are given as in Theorem \ref{theo:structure}. We then  have  

\begin{theorem}
Let $A$  be a singular irreducible M-matrix. Let $L_G$ and $U_G$ be the  unique matrices given by  \eqref{fspai1} and \eqref{fspai2} using $S_L$ and $S_U$ as given in \eqref{S_L-voll} and \eqref{S_U-voll}. Then the  matrix 
\begin{equation}
\hat A := U_GD^-L_G
\end{equation}
is a (1,2)-inverse of $A$
\end{theorem}
\begin{proof}
 With $I_F$ as in \eqref{eq:IF} we easily obtain
 \begin{equation*}
  \hat A A \hat A = U_GD^-L_GA U_GD^-L_G = U_GD^-I_FD^-L_G = U_GD^-L_G = \hat A.
 \end{equation*}
Next we consider $L_GA U_GD^-L_GAU_G$. We get  
 
 \begin{eqnarray*}
   L_GA U_GD^-L_GAU_G  = \left[ \begin{array}{ccccc}
                           d_1 & & & & \\
                           & \ddots & & & \\
                           & & d_{n-2} & & \\
                           & &  & d_{n-1} & b\\
                           & & & c & d_n
                          \end{array} \right]  \\
                          *   \left[ \begin{array}{ccccc}
                           d_1^{-1} & & & & \\
                           & \ddots & & & \\
                           & & d_{n-2}^{-1} & & \\
                           & &  & d_{n-1}^{-1} & \\
                           & & &  & 0
                          \end{array} \right] * 
                          \left[ \begin{array}{ccccc}
                           d_1 & & & & \\
                           & \ddots & & & \\
                           & & d_{n-2} & & \\
                           & &  & d_{n-1} & b\\
                           & & & c & d_n
                          \end{array} \right] \\
                          = 
                          \left[ \begin{array}{ccccc}
                           d_1 & & & & \\
                           & \ddots & & & \\
                           & & d_{n-2} & & \\
                           & &  & d_{n-1} & b\\
                           & & & c & d_n
                          \end{array} \right],
\end{eqnarray*}
 since $d_{n-1}d_n = cb$. Thus
 \begin{equation*}
  L_GA U_GD^-L_GAU_G = L_GA U_G.  
 \end{equation*}
But $L_G$ and $U_G$ are non-singular by Theorem \ref{theo:exis}, so we obtain

  \begin{equation*}
 A \hat  A A = A U_GD^-L_GA = A.  
 \end{equation*}
 Hence, $U_GD^-L_G$ is (1,2)-inverse of $A$. 
\end{proof}

But in general $U_GD^-L_G$ is not a Moore-Penrose inverse  of $A$, which can be seen by the following  example. 

\begin{example}
Let  $A$ be given by
\begin{equation*}
 A =  \left[ \begin{array}{ccc}
                           0.6667 & -0.3333 & -0.3333 \\
                           -0.2500 & 0.5000 & -0.2500 \\
                           -0.4000 & -0.4000 & 0.8000 
                          \end{array} \right] .
\end{equation*}
Then  we get 
\begin{eqnarray*}
L_G =  \left[ \begin{array}{ccc}
                           1.5000 & 0 & 0 \\
                           1.0000  & 2.6667 & 0 \\
                           1.0000 &  0  & 1.6667 
                          \end{array} \right] \quad 
U_G =  \left[ \begin{array}{ccc}
                           1.5000 & 1.3333 & 0.8333 \\
                           0 & 2.6667 & 0 \\
                           0 & 0  & 1.6667 
                          \end{array} \right] , \\
                          D^- =  \left[ \begin{array}{ccc}
                           0.6667 & 0 & 0 \\
                           0 & 0.3750 & 0 \\
                           0 & 0  & 0 
                          \end{array} \right] .
\end{eqnarray*}                        
Thus 
\begin{equation*}
 \hat A = U_GD^-L_G  = \left[ \begin{array}{ccc}
                           2.0000  & 1.3333 & 0 \\
                           1.0000 & 2.6667 & 0 \\
                           0 & 0  & 0 
                          \end{array} \right] .
\end{equation*}
But $ A \hat A \neq   (A \hat A)^T$  and   $ \hat A A \neq   (\hat A A)^T$, since 

\begin{equation*}
 \hat A A  = \left[ \begin{array}{ccc}
                           1  & 0 & -1 \\
                           0 & 1 & -1 \\
                           0 & 0  & 0 
                          \end{array} \right] \quad \text{and} \quad
   A \hat A = \left[ \begin{array}{ccc}
                           1  & 0 & 0 \\
                           0 & 1 & 0 \\
                           -1.2 & -1.6  & 0 
                          \end{array} \right]  .                      
                          \end{equation*}

 So we have constructed a   (1,2)-inverse, but not a (1,2,3,4)-inverse.
\end{example}

\section{Numerical examples}

In this section we  illustrate  the theoretical results  established in the previous sections. We considered four matrices obtained from Markov chain modeling. These  matrices were  also studied in \cite{BenT02} but of larger sizes, since there the parallel structure of sparse approximate inverses was considered. \\

\begin{center}
\begin{tabular}{|c||c|c|c|c|}\hline
matrix & size  n&\# nonzero elements & symmetric \cr\hline \hline
2D&121&441&no \cr
ncd&286&1606& no\cr
leaky&530&4186& no\cr
telecom&666&3091& no\cr
\hline
\end{tabular}\\
\end{center}
\begin{center}
Table 1, \\
Properties of the test matrices\\ 
\end{center}
\vspace{0.3cm}


In all   cases   we used the Bi-CGSTAB \cite{VdV92} method as a solver. The right hand  side is $b = 0$ and we start with a randomly chosen  vector $x_0$.
We iterate until the error tolerance  is below  $10^{-11}$ but stopped after 500 iterations. \\

We considered four preconditioner, the first one is $I$, i.e.  no
preconditioning. All  the others are FSAI preconditioner with  the following choice of the nonzero pattern:

\begin{itemize}
\item only the diagonal entries - fsaiD,
\item the same  non-zero pattern as the related lower (upper) part of $A$ - fsaiN,
\item band structure with bandwidth 5  in the left  and right preconditioner - fsaiB.
  \end{itemize}

The left and right preconditioner are constructed as above.


The first example is the 2D matrix  which comes from  a two-dimensional Markov chain model. 

\begin{center}
\begin{tabular}{|l||c|c|c|c|}\hline
&error&iterations&convergence \cr \hline \hline
$I$&$4.7493\cdot 10^{-12}$&73&yes \cr \hline 
fsaiD&$3.2091\cdot 10^{-13}$&31&yes \cr \hline
fsaiN&$2.3995\cdot 10^{-12}$&23&yes \cr \hline
fsaiB&$8.1504\cdot 10^{-12}$&25&yes \cr \hline
\end{tabular} \\
\end{center}
\begin{center}
Table 2 - 2D\\ 
\end{center}


This simple example shows already the potential of the FSAI preconditioning. The iteration number is one third compared to the unpreconditioined case.   

The next  example comes from a multiplexing model of a leaky bucket. \\

\begin{center}
\begin{tabular}{|l||c|c|c|c|}\hline
&error&iterations&convergence \cr \hline \hline
$I$&$1.2104\cdot 10^{-11}$&116&yes  \cr \hline 
fsaiD&$7.9897\cdot 10^{-12}$&109& yes \cr \hline
fsaiN&$2.7913\cdot 10^{-12}$&65& yes \cr \hline
fsaiB&$5.0377\cdot 10^{-12}$&111& yes \cr \hline
\end{tabular} \\
\end{center}
\begin{center}
Table 3 - leaky\\
\end{center}

An NCD queuing network  is the next example. \\

\begin{center}
\begin{tabular}{|l||c|c|c|c|}\hline
&error&iterations&convergence \cr \hline \hline
$I$&$9.7499\cdot 10^{-5}$&$>500$&no \cr \hline 
fsaiD&$1.3520\cdot 10^{-9}$&$>500$&no \cr \hline
fsaiN&$2.3031\cdot 10^{-12}$&84&yes \cr \hline
fsaiB&$7.2176\cdot 10^{-12}$&298&yes \cr \hline
\end{tabular} \\
\end{center}
\begin{center}
Table 4 - NCD\\
\end{center}


Finally,  a telecommunication model is considered. \\

\begin{center}
\begin{tabular}{|l||c|c|c|c|}\hline
&error&iterations&convergence \cr \hline \hline
$I$&$0.0159$&$>500$&no \cr \hline 
fsaiD&$8.9920\cdot 10^{-12}$&197&yes  \cr \hline
fsaiN&$6.6865\cdot 10^{-12}$&98& yes \cr \hline
fsaiB&$1.7722\cdot 10^{-12}$&194& yes \cr \hline
\end{tabular} \\
\vspace{0.3cm}
Table 5 - telecom
\end{center}

Especially,  the last two examples show  that the FSAI preconditioner with the same nonzero  structure as $A$ work very well  and seems  to be the best  choice. Moreover, it can be seen  that  just a diagonal preconditioning does not work at all. 


\printbibliography

@article{BenNFS01,
  title={Algebraic theory of multiplicative Schwarz methods},
  author={Benzi, Michele and Frommer, Andreas and Nabben, Reinhard and Szyld, Daniel B},
  journal={Numerische Mathematik},
  volume={89},
  number={4},
  pages={605--639},
  year={2001},
  publisher={Springer}
}

@article{ElsFNSS03,
title = {Conditions for strict inequality in comparisons of spectral radii of splittings of different matrices},
journal = {Linear Algebra and its Applications},
volume = {363},
pages = {65-80},
year = {2003},
note = {Special Issue on Nonnegative matrices, M-matrices and their generalizations},
issn = {0024-3795},
doi = {https://doi.org/10.1016/S0024-3795(01)00535-3},
url = {https://www.sciencedirect.com/science/article/pii/S0024379501005353},
author = {Ludwig Elsner and Andreas Frommer and Reinhard Nabben and Hans Schneider and Daniel B. Szyld}
}

@article{FerJP14,
  title={A generalized Block FSAI preconditioner for nonsymmetric linear systems},
  author={Ferronato, Massimiliano and Janna, Carlo and Pini, Giorgio},
  journal={Journal of Computational and Applied Mathematics},
  volume={256},
  pages={230--241},
  year={2014},
  publisher={Elsevier}
}

@article{KolY95,
  title={Factorized sparse approximate inverse preconditioning {II}: Solution of 3D FE systems on massively parallel computers},
  author={Kolotilina, L Yu and Yeremin, A Yu},
  journal={International Journal of High Speed Computing},
  volume={7},
  number={02},
  pages={191--215},
  year={1995},
  publisher={World Scientific}
}

@article{YerK00,
  title={Factorized sparse approximate inverse preconditionings. {III}. iterative construction of preconditioners},
  author={Yeremin, Alex Yu and Kolotilina, Lily Yu and Nikishin, AA},
  journal={Journal of Mathematical Sciences},
  volume={101},
  number={4},
  pages={3237--3254},
  year={2000},
  publisher={Springer}
}

@article{KolY99,
  title={Factorized sparse approximate inverse preconditionings. {IV}: Simple approaches to rising efficiency},
  author={Kolotilina, L Yu and Nikishin, Andy A and Yeremin, A Yu},
  journal={Numerical Linear Algebra with Applications},
  volume={6},
  number={7},
  pages={515--531},
  year={1999},
  publisher={Wiley Online Library}
}

@article{JanFG10,
  title={A block FSAI-ILU parallel preconditioner for symmetric positive definite linear systems},
  author={Janna, Carlo and Ferronato, Massimilano and Gambolati, Giuseppe},
  journal={SIAM Journal on Scientific Computing},
  volume={32},
  number={5},
  pages={2468--2484},
  year={2010},
  publisher={SIAM}
}

@article{JanF11,
  title={Adaptive pattern research for block FSAI preconditioning},
  author={Janna, Carlo and Ferronato, Massimiliano},
  journal={SIAM Journal on Scientific Computing},
  volume={33},
  number={6},
  pages={3357--3380},
  year={2011},
  publisher={SIAM}
}

@article{LivB12,
  title={Lean algebraic multigrid ({LAMG}): Fast graph {L}aplacian linear solver},
  author={Livne, Oren E and Brandt, Achi},
  journal={SIAM Journal on Scientific Computing},
  volume={34},
  number={4},
  pages={B499--B522},
  year={2012},
  publisher={SIAM}
}

@article{NapN16,
  title={An efficient multigrid method for graph {L}aplacian systems},
  author={Napov, Artem and Notay, Yvan},
  journal={Electron. Trans. Numer. Anal},
  volume={45},
  pages={201},
  year={2016}
}

@article{VdV92,
  title={{B}i-{CGSTAB}: A fast and smoothly converging variant of {B}i-{CG} for the solution of nonsymmetric linear systems},
  author={Van der Vorst, Henk A},
  journal={SIAM Journal on scientific and Statistical Computing},
  volume={13},
  number={2},
  pages={631--644},
  year={1992},
  publisher={SIAM}
}

@article{BolFFHK11,
  title={Algebraic multigrid methods for {L}aplacians of graphs},
  author={Bolten, Matthias and Friedhoff, Stephanie and Frommer, Andreas and Heming, Matthias and Kahl, Karsten},
  journal={Linear Algebra and its Applications},
  volume={434},
  number={11},
  pages={2225--2243},
  year={2011},
  publisher={Elsevier}
}

@article {HucS10,
    AUTHOR = {Huckle, T. and Sedlacek, M.},
     TITLE = {Smoothing and regularization with modified sparse approximate inverses},
   JOURNAL = {Journal of Electrical and Computer Engineering},
  FJOURNAL = {},
    VOLUME = {1},
      YEAR = {2010},
    NUMBER = {},
     PAGES = {930218},   
     ISSN = {},
   MRCLASS = {},
  MRNUMBER = {},
      }

@article{BroGMR01,
  title={Robust parallel smoothing for multigrid via sparse approximate inverses},
  author={Br{\"o}ker, Oliver and Grote, Marcus J and Mayer, Carsten and Reusken, Arnold},
  journal={SIAM Journal on Scientific Computing},
  volume={23},
  number={4},
  pages={1396--1417},
  year={2001},
  publisher={SIAM}
}

@article{BarBS99,
  title={An {MPI} implementation of the {SPAI} preconditioner on the {T3E}},
  author={Barnard, Stephen T and Bernardo, Luis M and Simon, Horst D},
  journal={The International Journal of High Performance Computing Applications},
  volume={13},
  number={2},
  pages={107--123},
  year={1999},
  publisher={Sage Publications Sage CA: Thousand Oaks, CA}
}

@article {Huc03,
    AUTHOR = {Huckle, T.},
     TITLE = {Factorized sparse approximate inverses for preconditioning},
   JOURNAL = {The Journal of Supercomputing},
  FJOURNAL = {},
    VOLUME = {25},
      YEAR = {2003},
    NUMBER = {},
     PAGES = {109--117},   
     ISSN = {},
   MRCLASS = {},
  MRNUMBER = {},
      }

@article {Huc00,
    AUTHOR = {Huckle, T.},
     TITLE = {Factorized sparse approximate inverses for preconditioning and smoothing},
   JOURNAL = {Selcuk Journal of Applied Mathematics},
  FJOURNAL = {},
    VOLUME = {1},
      YEAR = {2000},
    NUMBER = {63},
     PAGES = {63--70},   
     ISSN = {},
   MRCLASS = {},
  MRNUMBER = {},
      }

@article {ErlN17,
    AUTHOR = {Erlangga, Y.A. and Nabben, R.},
     TITLE = {On the convergence of two‐level {K}rylov methods for singular symmetric systems},
   JOURNAL = {Numer. Linear. Algebra Appl.},
  FJOURNAL = {Numerical Linear Algebra with Applications},
    VOLUME = {24},
      YEAR = {2017},
    NUMBER = {6},
     PAGES = {e2108},   
     ISSN = {},
   MRCLASS = {},
  MRNUMBER = {},
      }

@article {LudNT16,
    AUTHOR = {Ludwig, E. and Nabben, R. and Tang, J.M.},
     TITLE = {Deflation and projection methods applied to symmetric positive semi-definite systems},
   JOURNAL = {Linear Algebra Appl.},
  FJOURNAL = {Linear Algebra and its Applications},
    VOLUME = {489},
      YEAR = {2016},
    NUMBER = {},
     PAGES = {253--273},   
      ISSN = {},
   MRCLASS = {},
  MRNUMBER = {},
      }

@article {BenT98,
    AUTHOR = {Benzi, Michele and T\r{u}ma, Miroslav},
     TITLE = {Numerical experiments with two approximate inverse
              preconditioners},
   JOURNAL = {BIT},
  FJOURNAL = {BIT. Numerical Mathematics},
    VOLUME = {38},
      YEAR = {1998},
    NUMBER = {2},
     PAGES = {234--241},
      ISSN = {0006-3835},
   MRCLASS = {65F35 (65F10 65F50)},
  MRNUMBER = {1638175},
       DOI = {10.1007/BF02512364},
       URL = {https://doi.org/10.1007/BF02512364},
}

@article {KolY93,
    AUTHOR = {Kolotilina, L. Yu. and Yeremin, A. Yu.},
     TITLE = {Factorized sparse approximate inverse preconditionings. {I}.
              {T}heory},
   JOURNAL = {SIAM J. Matrix Anal. Appl.},
  FJOURNAL = {SIAM Journal on Matrix Analysis and Applications},
    VOLUME = {14},
      YEAR = {1993},
    NUMBER = {1},
     PAGES = {45--58},
      ISSN = {0895-4798},
   MRCLASS = {65F50 (65F35 65Y05)},
  MRNUMBER = {1199543},
MRREVIEWER = {Dao\ Sheng\ Zheng},
       DOI = {10.1137/0614004},
       URL = {https://doi.org/10.1137/0614004},
}

@article {KhaKNY01,
    AUTHOR = {Kharchenko, S. A. and Kolotilina, L. Yu. and Nikishin, A. A.
              and Yeremin, A. Yu.},
     TITLE = {A robust {AINV}-type method for constructing sparse
              approximate inverse preconditioners in factored form},
   JOURNAL = {Numer. Linear Algebra Appl.},
  FJOURNAL = {Numerical Linear Algebra with Applications},
    VOLUME = {8},
      YEAR = {2001},
    NUMBER = {3},
     PAGES = {165--179},
      ISSN = {1070-5325,1099-1506},
   MRCLASS = {65F10},
  MRNUMBER = {1817794},
       DOI = {10.1002/1099-1506(200104/05)8:3<165::AID-NLA235>3.0.CO;2-9},
       URL =  {https://doi.org/10.1002/1099-1506(200104/05)8:3<165::AID-NLA235>3.0.CO;2-9},
}

@article {BenB03,
    AUTHOR = {Benzi, Michele and Bertaccini, Daniele},
     TITLE = {Approximate inverse preconditioning for shifted linear
              systems},
   JOURNAL = {BIT},
  FJOURNAL = {BIT. Numerical Mathematics},
    VOLUME = {43},
      YEAR = {2003},
    NUMBER = {2},
     PAGES = {231--244},
      ISSN = {0006-3835,1572-9125},
   MRCLASS = {65F10},
  MRNUMBER = {2010362},
       DOI = {10.1023/A:1026089811044},
       URL = {https://doi.org/10.1023/A:1026089811044},
}

@incollection {Mey82,
    AUTHOR = {Meyer, Jr., C. D.},
     TITLE = {Analysis of finite {M}arkov chains by group inversion
              techniques},
 BOOKTITLE = {Recent applications of generalized inverses},
    SERIES = {Res. Notes in Math.},
    VOLUME = {66},
     PAGES = {50--81},
 PUBLISHER = {Pitman, Boston, Mass.-London},
      YEAR = {1982},
      ISBN = {0-273-08550-6},
   MRCLASS = {60J10 (15A09)},
  MRNUMBER = {666723},
MRREVIEWER = {Adolf\ Rhodius},
}

@book {KirN13,
    AUTHOR = {Kirkland, Stephen J. and Neumann, Michael},
     TITLE = {Group inverses of {M}-matrices and their applications},
    SERIES = {Chapman \& Hall/CRC Applied Mathematics and Nonlinear Science
              Series},
 PUBLISHER = {CRC Press, Boca Raton, FL},
      YEAR = {2013},
     PAGES = {xvi+316},
      ISBN = {978-1-4398-8858-2},
   MRCLASS = {15-02 (15A09 15B48 60J10)},
  MRNUMBER = {3185162},
MRREVIEWER = {Jian\ Long\ Chen},
}

@article {GroH97,
    AUTHOR = {Grote, Marcus J. and Huckle, Thomas},
     TITLE = {Parallel preconditioning with sparse approximate inverses},
   JOURNAL = {SIAM J. Sci. Comput.},
  FJOURNAL = {SIAM Journal on Scientific Computing},
    VOLUME = {18},
      YEAR = {1997},
    NUMBER = {3},
     PAGES = {838--853},
      ISSN = {1064-8275},
   MRCLASS = {65Fxx (65Y05)},
  MRNUMBER = {1443644},
MRREVIEWER = {Zahari\ Zlatev},
       DOI = {10.1137/S1064827594276552},
       URL = {https://doi.org/10.1137/S1064827594276552},
}

@article {BolHK16,
    AUTHOR = {Bolten, Matthias and Huckle, Thomas K. and Kravvaritis,
              Christos D.},
     TITLE = {Sparse matrix approximations for multigrid methods},
   JOURNAL = {Linear Algebra Appl.},
  FJOURNAL = {Linear Algebra and its Applications},
    VOLUME = {502},
      YEAR = {2016},
     PAGES = {58--76},
      ISSN = {0024-3795,1873-1856},
   MRCLASS = {65N55 (65F10 65F15 65F50)},
  MRNUMBER = {3490785},
MRREVIEWER = {Svetozar\ D.\ Margenov},
       DOI = {10.1016/j.laa.2015.11.008},
       URL = {https://doi.org/10.1016/j.laa.2015.11.008},
}

@article {AnzHBD18,
    AUTHOR = {Anzt, Hartwig and Huckle, Thomas K. and Br\"ackle, J\"urgen
              and Dongarra, Jack},
     TITLE = {Incomplete sparse approximate inverses for parallel
              preconditioning},
   JOURNAL = {Parallel Comput.},
  FJOURNAL = {Parallel Computing. Systems \& Applications},
    VOLUME = {71},
      YEAR = {2018},
     PAGES = {1--22},
      ISSN = {0167-8191,1872-7336},
   MRCLASS = {65F08 (65Y05)},
  MRNUMBER = {3735016},
       DOI = {10.1016/j.parco.2017.10.003},
       URL = {https://doi.org/10.1016/j.parco.2017.10.003},
}

@article {BenT99,
    AUTHOR = {Benzi, Michele and T\r{u}ma, Miroslav},
    TITLE = {A comparative study of sparse approximate inverse
              preconditioners},
   JOURNAL = {Appl. Numer. Math.},
  FJOURNAL = {Applied Numerical Mathematics. An IMACS Journal},
    VOLUME = {30},
      YEAR = {1999},
    NUMBER = {2-3},
     PAGES = {305--340},
      ISSN = {0168-9274,1873-5460},
   MRCLASS = {65F50 (65F35)},
  MRNUMBER = {1688632},
       DOI = {10.1016/S0168-9274(98)00118-4},
       URL = {https://doi.org/10.1016/S0168-9274(98)00118-4},
}

@article {BenT98a,
    AUTHOR = {Benzi, Michele and T\r{u}ma, Miroslav},
     TITLE = {A sparse approximate inverse preconditioner for nonsymmetric
              linear systems},
   JOURNAL = {SIAM J. Sci. Comput.},
  FJOURNAL = {SIAM Journal on Scientific Computing},
    VOLUME = {19},
      YEAR = {1998},
    NUMBER = {3},
     PAGES = {968--994},
      ISSN = {1064-8275,1095-7197},
   MRCLASS = {65F35 (65F50)},
  MRNUMBER = {1616710},
       DOI = {10.1137/S1064827595294691},
       URL = {https://doi.org/10.1137/S1064827595294691},
}

@article {BenMT96,
    AUTHOR = {Benzi, Michele and Meyer, Carl D. and T\r{u}ma, Miroslav},
     TITLE = {A sparse approximate inverse preconditioner for the conjugate
              gradient method},
   JOURNAL = {SIAM J. Sci. Comput.},
  FJOURNAL = {SIAM Journal on Scientific Computing},
    VOLUME = {17},
      YEAR = {1996},
    NUMBER = {5},
     PAGES = {1135--1149},
      ISSN = {1064-8275},
   MRCLASS = {65F50 (65F10 65Y05)},
  MRNUMBER = {1404865},
       DOI = {10.1137/S1064827594271421},
       URL = {https://doi.org/10.1137/S1064827594271421},
}

@article {BenCT00,
    AUTHOR = {Benzi, Michele and Cullum, Jane K. and T\r{u}ma, Miroslav},
     TITLE = {Robust approximate inverse preconditioning for the conjugate
              gradient method},
   JOURNAL = {SIAM J. Sci. Comput.},
  FJOURNAL = {SIAM Journal on Scientific Computing},
    VOLUME = {22},
      YEAR = {2000},
    NUMBER = {4},
     PAGES = {1318--1332},
      ISSN = {1064-8275,1095-7197},
   MRCLASS = {65F10 (65F50 65N22)},
  MRNUMBER = {1787297},
MRREVIEWER = {Zahari\ Zlatev},
       DOI = {10.1137/S1064827599356900},
       URL = {https://doi.org/10.1137/S1064827599356900},
}

@article {BenT02,
    AUTHOR = {Benzi, Michele and T\r{u}ma, Miroslav},
     TITLE = {A parallel solver for large-scale {M}arkov chains},
   JOURNAL = {Appl. Numer. Math.},
  FJOURNAL = {Applied Numerical Mathematics. An IMACS Journal},
    VOLUME = {41},
      YEAR = {2002},
    NUMBER = {1},
     PAGES = {135--153},
      ISSN = {0168-9274,1873-5460},
   MRCLASS = {65F10 (65C40)},
  MRNUMBER = {1908754},
MRREVIEWER = {Jiu\ Ding},
       DOI = {10.1016/S0168-9274(01)00116-7},
       URL = {https://doi.org/10.1016/S0168-9274(01)00116-7},
}

@article {BenT03,
    AUTHOR = {Benzi, Michele and  T\r{u}ma, Miroslav},
     TITLE = {A robust incomplete factorization preconditioner for positive
              definite matrices},
      NOTE = {Preconditioning, 2001 (Tahoe City, CA)},
   JOURNAL = {Numer. Linear Algebra Appl.},
  FJOURNAL = {Numerical Linear Algebra with Applications},
    VOLUME = {10},
      YEAR = {2003},
    NUMBER = {5-6},
     PAGES = {385--400},
      ISSN = {1070-5325,1099-1506},
   MRCLASS = {65F10},
  MRNUMBER = {2008366},
       DOI = {10.1002/nla.320},
       URL = {https://doi.org/10.1002/nla.320},
}

@book {CamM09,
    AUTHOR = {Campbell, Stephen L. and Meyer, Carl D.},
     TITLE = {Generalized inverses of linear transformations},
    SERIES = {Classics in Applied Mathematics},
    VOLUME = {56},
      NOTE = {},
 PUBLISHER = {Society for Industrial and Applied Mathematics (SIAM),
              Philadelphia, PA},
      YEAR = {2009},
     PAGES = {xx+272},
      ISBN = {978-0-898716-71-9},
   MRCLASS = {15A09 (01A75 15-02 47A05 65F05 65F20)},
  MRNUMBER = {3396208},
       DOI = {10.1137/1.9780898719048.ch0},
}

@book {BerP79,
    AUTHOR = {Berman, Abraham and Plemmons, Robert J.},
     TITLE = {Nonnegative matrices in the mathematical sciences},
    SERIES = {Computer Science and Applied Mathematics},
 PUBLISHER = {Academic Press [Harcourt Brace Jovanovich, Publishers], New
              York-London},
      YEAR = {1979},
     PAGES = {xviii+316},
      ISBN = {0-12-092250-9},
   MRCLASS = {15A48 (15-02 60J05 90C33)},
  MRNUMBER = {544666},
MRREVIEWER = {R.\ S.\ Varga},
}

@article {FieP62,
    AUTHOR = {Fiedler, Miroslav and Pt\'ak, Vlastimil},
     TITLE = {On matrices with non-positive off-diagonal elements and
              positive principal minors},
   JOURNAL = {Czechoslovak Math. J.},
  FJOURNAL = {Czechoslovak Mathematical Journal},
    VOLUME = {12(87)},
      YEAR = {1962},
     PAGES = {382--400},
      ISSN = {0011-4642},
   MRCLASS = {15.60},
  MRNUMBER = {142565},
MRREVIEWER = {Ky\ Fan},
}

@article {NabV95,
    AUTHOR = {Nabben, Reinhard and Varga, Richard S.},
     TITLE = {On classes of inverse {$Z$}-matrices},
      NOTE = {Special issue honoring Miroslav Fiedler and Vlastimil Pt\'ak},
   JOURNAL = {Linear Algebra Appl.},
  FJOURNAL = {Linear Algebra and its Applications},
    VOLUME = {223/224},
      YEAR = {1995},
     PAGES = {521--552},
      ISSN = {0024-3795,1873-1856},
   MRCLASS = {15A09 (15A48)},
  MRNUMBER = {1340709},
MRREVIEWER = {Ronald\ L.\ Smith},
       DOI = {10.1016/0024-3795(95)92718-8},
       URL = {https://doi.org/10.1016/0024-3795(95)92718-8},
}

@article {Nab97,
    AUTHOR = {Nabben, Reinhard},
     TITLE = {{$Z$}-matrices and inverse {$Z$}-matrices},
   JOURNAL = {Linear Algebra Appl.},
  FJOURNAL = {Linear Algebra and its Applications},
    VOLUME = {256},
      YEAR = {1997},
     PAGES = {31--48},
      ISSN = {0024-3795,1873-1856},
   MRCLASS = {15A57},
  MRNUMBER = {1438224},
MRREVIEWER = {Thomas\ L.\ Markham},
       DOI = {10.1016/S0024-3795(97)81111-1},
       URL = {https://doi.org/10.1016/S0024-3795(97)81111-1},
}

@article {FieM92,
    AUTHOR = {Fiedler, Miroslav and Markham, Thomas L.},
     TITLE = {A classification of matrices of class {$Z$}},
   JOURNAL = {Linear Algebra Appl.},
  FJOURNAL = {Linear Algebra and its Applications},
    VOLUME = {173},
      YEAR = {1992},
     PAGES = {115--124},
      ISSN = {0024-3795,1873-1856},
   MRCLASS = {15A57},
  MRNUMBER = {1170507},
       DOI = {10.1016/0024-3795(92)90425-A},
       URL = {https://doi.org/10.1016/0024-3795(92)90425-A},
}

@article {FroNS08,
    AUTHOR = {Frommer, Andreas and Nabben, Reinhard and Szyld, Daniel B.},
     TITLE = {Convergence of stationary iterative methods for {H}ermitian
              semidefinite linear systems and applications to {S}chwarz
              methods},
   JOURNAL = {SIAM J. Matrix Anal. Appl.},
  FJOURNAL = {SIAM Journal on Matrix Analysis and Applications},
    VOLUME = {30},
      YEAR = {2008},
    NUMBER = {2},
     PAGES = {925--938},
      ISSN = {0895-4798,1095-7162},
   MRCLASS = {65F10 (65F20)},
  MRNUMBER = {2443978},
MRREVIEWER = {Giuseppe\ Rodriguez},
       DOI = {10.1137/080714038},
       URL = {https://doi.org/10.1137/080714038},
}

@article {NABS07,
    AUTHOR = {Nabben, Reinhard and Szyld, Daniel B.},
     TITLE = {Schwarz iterations for symmetric positive semidefinite
              problems},
   JOURNAL = {SIAM J. Matrix Anal. Appl.},
  FJOURNAL = {SIAM Journal on Matrix Analysis and Applications},
    VOLUME = {29},
      YEAR = {2006/07},
    NUMBER = {1},
     PAGES = {98--116},
      ISSN = {0895-4798,1095-7162},
   MRCLASS = {65F10},
  MRNUMBER = {2288016},
       DOI = {10.1137/050644203},
       URL = {https://doi.org/10.1137/050644203},
}

@article {BruPS05,
    AUTHOR = {Bru, Rafael and Pedroche, Francisco and Szyld, Daniel B.},
     TITLE = {Additive {S}chwarz iterations for {M}arkov chains},
   JOURNAL = {SIAM J. Matrix Anal. Appl.},
  FJOURNAL = {SIAM Journal on Matrix Analysis and Applications},
    VOLUME = {27},
      YEAR = {2005},
    NUMBER = {2},
     PAGES = {445--458},
      ISSN = {0895-4798,1095-7162},
   MRCLASS = {65F10},
  MRNUMBER = {2179682},
MRREVIEWER = {Vladimir\ B.\ Larin},
       DOI = {10.1137/040616541},
       URL = {https://doi.org/10.1137/040616541},
}

@article {MarS93,
    AUTHOR = {Marek, Ivo and Szyld, Daniel B.},
     TITLE = {Iterative and semi-iterative methods for computing stationary
              probability vectors of {M}arkov operators},
   JOURNAL = {Math. Comp.},
  FJOURNAL = {Mathematics of Computation},
    VOLUME = {61},
      YEAR = {1993},
    NUMBER = {204},
     PAGES = {719--731},
      ISSN = {0025-5718,1088-6842},
   MRCLASS = {65J10 (15A48 47A50 47B65 60J10)},
  MRNUMBER = {1192973},
MRREVIEWER = {I.\ K.\ Marek},
       DOI = {10.2307/2153249},
       URL = {https://doi.org/10.2307/2153249},
}

@article {Szy94,
    AUTHOR = {Szyld, Daniel B.},
     TITLE = {Equivalence of conditions for convergence of iterative methods
              for singular equations},
   JOURNAL = {Numer. Linear Algebra Appl.},
  FJOURNAL = {Numerical Linear Algebra with Applications},
    VOLUME = {1},
      YEAR = {1994},
    NUMBER = {2},
     PAGES = {151--154},
      ISSN = {1070-5325,1099-1506},
   MRCLASS = {65J10 (47A50)},
  MRNUMBER = {1277799},
       DOI = {10.1002/nla.1680010206},
       URL = {https://doi.org/10.1002/nla.1680010206},
}

@article {MarS00,
    AUTHOR = {Marek, Ivo and Szyld, Daniel B.},
     TITLE = {Comparison theorems for the convergence factor of iterative
              methods for singular matrices},
      NOTE = {Conference Celebrating the 60th Birthday of Robert J. Plemmons
              (Winston-Salem, NC, 1999)},
   JOURNAL = {Linear Algebra Appl.},
  FJOURNAL = {Linear Algebra and its Applications},
    VOLUME = {316},
YEAR = {2000},}

@article {Ben02,
    AUTHOR = {Benzi, Michele},
     TITLE = {Preconditioning techniques for large linear systems: a survey},
   JOURNAL = {J. Comput. Phys.},
  FJOURNAL = {Journal of Computational Physics},
    VOLUME = {182},
      YEAR = {2002},
    NUMBER = {2},
     PAGES = {418--477},
      ISSN = {0021-9991},
   MRCLASS = {65F10 (65F50)},
  MRNUMBER = {1941848},
       DOI = {10.1006/jcph.2002.7176},
       URL = {https://doi.org/10.1006/jcph.2002.7176},
}

@article {MenN08a,
    AUTHOR = {Mense, Christian and Nabben, Reinhard},
     TITLE = {On algebraic multilevel methods for non-symmetric
              systems---convergence results},
   JOURNAL = {Electron. Trans. Numer. Anal.},
  FJOURNAL = {Electronic Transactions on Numerical Analysis},
    VOLUME = {30},
      YEAR = {2008},
     PAGES = {323--345},
   MRCLASS = {65F10},
  MRNUMBER = {2487765},
}

@article {MenN08b,
    AUTHOR = {Mense, C. and Nabben, R.},
     TITLE = {On algebraic multi-level methods for non-symmetric
              systems---comparison results},
   JOURNAL = {Linear Algebra Appl.},
  FJOURNAL = {Linear Algebra and its Applications},
    VOLUME = {429},
      YEAR = {2008},
    NUMBER = {10},
     PAGES = {2567--2588},
      ISSN = {0024-3795},
   MRCLASS = {65N55},
  MRNUMBER = {2456796},
MRREVIEWER = {Eveline Rosseel},
       DOI = {10.1016/j.laa.2008.04.045},
       URL = {https://doi.org/10.1016/j.laa.2008.04.045},
}

@article {KehN16,
    AUTHOR = {Kehl, Ren\'{e} and Nabben, Reinhard},
     TITLE = {Avoiding singular coarse grid systems},
   JOURNAL = {Linear Algebra Appl.},
  FJOURNAL = {Linear Algebra and its Applications},
    VOLUME = {507},
      YEAR = {2016},
     PAGES = {137--152},
      ISSN = {0024-3795,1873-1856},
   MRCLASS = {65F10 (65F50 65N22 65N55)},
  MRNUMBER = {3536948},
MRREVIEWER = {Bruno\ Carpentieri},
       DOI = {10.1016/j.laa.2016.05.025},
       URL = {https://doi.org/10.1016/j.laa.2016.05.025},
}

\end{document}